\documentclass{article}
\usepackage{graphicx} 
\usepackage[utf8]{inputenc}
\usepackage{titlesec}
\titleformat{\section}
{\normalfont\scshape}
{\thesection.}{0.5em}{\centering}
\titleformat{\subsection}
{\small\scshape}
{\thesubsection.}{0.5em}{\centering}
\titleformat{\subsubsection}
{\small\scshape}
{\thesubsubsection.}{0.5em}{\centering}
\usepackage{enumerate}   
\usepackage{tikz-cd}
\usepackage{tikz}
\usepackage{comment}
\usepackage{amsfonts}
\usepackage{indentfirst}
\renewcommand{\abstract}{\small{\section*{\abstractname}}}
\usepackage{amsmath}
\usepackage{amssymb}
\usepackage{amsthm}
\usepackage{mathrsfs}
\newtheorem{theorem}{Theorem}
\newtheorem*{theorem*}{Theorem}
\newtheorem{lemma}{Lemma}

\newtheorem*{corollary*}{Corollary}
\newtheorem{definition}{Definition}

\theoremstyle{remark}
\newtheorem{remark}{Remark}
\newtheorem{example}{Example}
\usepackage{stmaryrd}
\newcommand{\hooklongrightarrow}{\lhook\joinrel\longrightarrow}

\usepackage[
backend=biber,
style=alphabetic,
sorting=ynt
]{biblatex}
\title{\Large{\textbf{Limit theorems for $\sigma$-localized Émery convergence}}}
\author{\normalsize{\textsc{Vasily Melnikov}}}
\date{\normalsize{\textsc{October 2024}}}

\begin{document}

\maketitle
\begin{abstract}
     Given a bounded sequence $\{X^{n}\}_{n}$ of semimartingales on a time interval $[0,T]$, we find a sequence of convex combinations $\{Y^{n}\}_{n}$ and a limiting semimartingale $Y$ such that $\{Y^{n}\}_{n}$ converges to $Y$ in a $\sigma$-localized modification of the Émery topology. More precisely, $\{Y^{n}\}_{n}$ converges to $Y$ in the Émery topology on an increasing sequence $\{D_{n}\}_{n}$ of predictable sets covering $\Omega\times[0,T]$. We also prove some technical variants of this theorem, including a version where the complement of $\{D_{n}\}_{n}$ forms a disjoint sequence. Applications include a complete characterization of sequences admitting convex combinations converging in the Émery topology, and a supermartingale counterpart of Helly's selection theorem.
\end{abstract}
\section{Introduction}\label{sec:intro}
Extracting a limit from a sequence $\{x_{n}\}_{n}$ in a topological space $X$ is a foundational procedure in mathematical analysis. When $X$ is given a linear structure over $\mathbb{R}$, it is often just as useful to find convergent convex combinations $y_{m}\in\mathrm{co}\{x_{n}:n\geq m\}$, and characterizing the well-posedness of this weaker problem is therefore of significant interest. Indeed, when $X$ is the space of random variables equipped with the topology of convergence in probability, convex compactness theorems bridge arbitrage theory and the Bichteler-Dellacherie characterization of semimartingales (see \cite{arb}), provide a simple proof of the Doob-Meyer decomposition (see \cite{doobmey}), and are a key tool in portfolio optimization (see \cite{karatzit,kramsch}).
\par
Suitable generalizations to classes of stochastic processes have been explored for finite variation predictable processes (see \cite{schcamp}), nonnegative martingales (see \cite{martconvcomp}), and Hardy martingales (see \cite{del-sch,melnikov}). However, in this context, convergence is understood in fairly weak senses; for example, convex compactness results are formulated using the notion of Fatou convergence of stochastic processes introduced by Föllmer and Kramkov \cite{optdecompconstrain}, or by considering pointwise convergence on $\Omega\times[0,\infty)$ (modulo negligible sets).
\par
Many applications require stronger notions of convergence, such as convergence in the semimartingale topology introduced by Émery \cite{emery}. Indeed, semimartingale convergence has applications to portfolio optimization (see \cite{kardem}), superreplication under constraints (see \cite{constraintappl,constraintapplv2}), and arbitrage theory (see \cite{ftap,emeryconvcomp}).
\par
To this end, we establish a Bolzano-Weierstrass-type principle for the semimartingale topology, allowing one to pass from a bounded sequence of semimartingales to a semimartingale limit.
\begin{theorem*}
    Let $\{X^{n}\}_{n}$ be a sequence of semimartingales such that
    \begin{equation}\label{eq:boundednesscond}
        \mathrm{co}\left\{\vert{(\xi\cdot X^{n})_{T}\vert}:n\in\mathbb{N},\xi\textrm{ is predictable and }\vert{\xi}\vert\leq1\right\}
    \end{equation}
    is bounded in probability.  Then there exists an increasing sequence $\{D_{n}\}_{n}$ of predictable sets increasing to $\Omega\times[0,T]$, a semimartingale $Y$, and $Y^{n}\in\mathrm{co}\{X^{m}:m\geq n\}$, such that $\{(\mathbf{1}_{D_{m}}\cdot Y^{n})\}_{n}$ converges to $(\mathbf{1}_{D_{m}}\cdot Y)$ in the semimartingale topology for each $m\in\mathbb{N}$.
\end{theorem*}
A sequence $\{X^{n}\}_{n}$ satisfying boundedness of (\ref{eq:boundednesscond}) may fail to admit $Y^{n}\in\mathrm{co}\{X^{m}:m\geq n\}$ converging in the semimartingale topology. The notion of convergence we use is therefore adapted to deal with such counterexamples.
\par
In analogy with a seminal paper of Kallsen \cite{sigmaloc}, we term this notion the $\sigma$-localization of Émery convergence. Given a class $\mathscr{C}$ of semimartingales, $X$ belongs to the $\sigma$-localization of $\mathscr{C}$ if there exists an increasing sequence $\{D_{n}\}_{n}$ of predictable sets increasing to $\Omega\times[0,T]$ such that $(\mathbf{1}_{D_{n}}\cdot X)\in\mathscr{C}$ for each $n$. Immediate are the parallels between our result, and $\sigma$-localization.
\par
$\sigma$-localization is natural in this context, since many non-$\sigma$-localized concepts in stochastic analysis and mathematical finance have proved inadequately broad in full generality (especially in the presence of unbounded jumps). For example, the fundamental theorem of asset pricing without technical boundedness assumptions \textit{must} make use of $\sigma$-martingales (see Example 2.3, \cite{unboundftap}), as opposed to the more convenient notion of a local martingale. Likewise, necessity dictates displacing Émery convergence with its $\sigma$-localization.
\par
The strength of using $\sigma$-localized Émery convergence is that it ensures the limiting process $Y$ obtained from $Y^{n}\in\mathrm{co}\{X^{m}:m\geq n\}$ is always a semimartingale and provides a sizeable class of integrands $\xi$ for which $\{(\xi\cdot {Y^{n}})_{T}\}_{n}$ converges to $(\xi\cdot Y)_{T}$. These properties are not shared by many modes of convergence considered in the literature (including those suggested by \cite{martconvcomp,optdecompconstrain,schcamp}), and especially so when Bolzano-Weierstrass-type principles hold. Admittedly, specific classes of semimartingales—such as the $L^{2}$-martingales—can have remarkable compactness properties, even under the Émery topology. However, our goal is to see how much compactness one can salvage from the pathology of general semimartingales, so consideration of the ideal situation is not helpful.
\par
Our results are parallel to classical measure-theoretic tools such as the Kadec-Pełczyński decomposition \cite{ogkp}. We develop this correspondence further, showing how a crucial corollary to the Kadec-Pełczyński decomposition, which completely characterizes when one can pass to $L^{1}$-convergent convex combinations, has an analogue for semimartingales. It gives the precise necessary and sufficient conditions for a sequence of semimartingales to admit convex combinations converging in the Émery topology.
\par
Implicit in all of the above results is a connatural boundedness condition, asking that the set
\begin{equation*}
        \mathrm{co}\left\{\vert{(\xi\cdot X^{n})_{T}\vert}:n\in\mathbb{N},\xi\textrm{ is predictable and }\vert{\xi}\vert\leq1\right\}
\end{equation*}
is bounded in probability. For our main results to have wide applicability, it is crucial to determine when their preconditions are satisfied, and we therefore also investigate the above boundedness condition. Applications include a version of Helly's selection theorem for supermartingales, which refines or generalizes results due to \cite{schcamp,martconvcomp,melnikov}.
\par
The paper is structured as follows. Section \ref{sec:prelim} reviews some preliminaries and establishes our notation. Section \ref{sec:main} states our main results, whose proofs are in Section \ref{sec:proofs}. Section \ref{sec:emeryconvg} contains an application of our results from Section \ref{sec:main} to a semimartingale counterpart of the Dunford-Pettis theorem. Section \ref{sec:super} investigates the boundedness conditions implicit in our main results, which is then applied to deduce a version of Helly's selection theorem for supermartingales. Appendix \ref{app:counter} probes some limitations of our main results via two counterexamples.
\section{Preliminaries}\label{sec:prelim}
Let $(\Omega,\mathscr{F},\mathbb{P})$ be a probability space. Fix a terminal date $T\in(0,\infty)$, and let $\mathbb{F}=\{\mathscr{F}_{t}:t\in[0,T]\}$ be a filtration of sub-$\sigma$-algebras of $\mathscr{F}$ such that the following modifications of the usual conditions hold:
\begin{enumerate}
    \item $\mathscr{F}_{0}$ is the $\mathbb{P}$-completion of the trivial $\sigma$-algebra $\{\emptyset,\Omega\}$.
    \item $\mathbb{F}$ is right continuous, i.e.,
\begin{equation*}
    \bigcap_{s>t}\mathscr{F}_{s}=\mathscr{F}_{t},
\end{equation*}
for each $t\in[0,T)$.
\end{enumerate}
Given a measure $\mathbb{Q}$ on $\mathscr{F}$ we write $\mathbb{Q}\ll\mathbb{P}$ if $\mathbb{P}(E)=0$ implies that $\mathbb{Q}(E)=0$ for $E\in\mathscr{F}$; $\mathbb{P}\ll\mathbb{Q}$ is defined similarly. If $\mathbb{P}\ll\mathbb{Q}$ and $\mathbb{Q}\ll\mathbb{P}$ we say that $\mathbb{P}$ and $\mathbb{Q}$ are \textit{equivalent}, and write $\mathbb{P}\sim\mathbb{Q}$.
\par
A set $E\subset\Omega\times[0,T]$ is said to be \textit{evanescent} if
\begin{equation*}
    \{\omega\in\Omega:\exists t\in[0,T]\textrm{ such that }(t,\omega)\in E\}
\end{equation*}
is a $\mathbb{P}$-null set. Stochastic processes are viewed modulo evanescent sets. A \textit{càdlàg} function is a function $\xi:[0,T]\longrightarrow\mathbb{R}$ which is right continuous, and has left limits at each point in its domain. We generally assume from this point onwards that all stochastic processes are càdlàg (at least up to evanescence).
\par
A \textit{semimartingale} is an adapted process $X$ which decomposes as $X=M+A$, where $M$ is a càdlàg local martingale under $\mathbb{P}$, and $A$ is an adapted càdlàg finite variation process. By the Bichteller-Dellacherie theorem (see, for example, \cite{arb}), the semimartingales are essentially the largest class of processes for which stochastic integration can be reasonably defined. Given semimartingales $X$ and $Y$, their quadratic covariation is denoted by $[X,Y]$.
\par
Let $\mathscr{P}$ denote the predictable $\sigma$-algebra on $\Omega\times[0,T]$. A process is \textit{predictable} if it is $\mathscr{P}$-measurable. A stopping time $\tau$ is \textit{predictable} if $\llbracket{\tau}\rrbracket\in\mathscr{P}$. Denote by $\mathscr{P}(1)$ the space of processes $\xi$ which are predictable, and satisfy $\vert{\xi}\vert\leq1$.
\par
Let $X$ be a semimartingale, and let $\xi$ be an $X$-integrable predictable process. Both $(\xi\cdot X)$ and $\int_{0}^{\cdot}\xi dX$ are used to denote the Itô stochastic integral of $\xi$ with respect to $X$.
\par
On the space of semimartingales, it is necessary to consider various modes of convergence. As such, we consider the \textit{semimartingale topology} (introduced by Émery \cite{emery}).
\begin{definition}
    The \textit{semimartingale topology} is the topology on the space of (equivalence classes modulo evanescence of) semimartingales induced by the metric $\mathbf{D}$ defined as
    \begin{equation*}
        \mathbf{D}(X,Y)=\vert{X_{0}-Y_{0}}\vert+\sup_{\xi\in\mathscr{P}(1)}\int_{\Omega}\left(\xi\cdot(X-Y)\right)^{\ast}_{T}\wedge1d\mathbb{P},
    \end{equation*}
    for any semimartingales $X$ and $Y$.
\end{definition}
Let $\mathscr{C}$ be a class of semimartingales. Following (Definition 2.1, \cite{sigmaloc}), we say that $X$ belongs to the \textit{$\sigma$-localization} of $\mathscr{C}$ if there exists an increasing sequence $\{D_{n}\}_{n}$ of predictable sets such that $(\mathbf{1}_{D_{n}}\cdot X)\in\mathscr{C}$ for each $n$, and $\bigcup_{n=1}^{\infty}D_{n}=\Omega\times[0,T]$. By taking $D_{n}=\llbracket{0,T_{n}}\rrbracket$ for some localizing sequence $\{T_{n}\}_{n}$ one recovers the notion of localization.
\par
Let $\mathscr{V}^{d}$ denote the set of finite variation pure-jump semimartingales, i.e., $X\in\mathscr{V}^{d}$ if, and only if,
\begin{equation*}
    X=X_{0}+x\ast\mu^{X},
\end{equation*}
where $\mu^{X}$ is the jump measure of $X$, and $\ast$ denotes integration with respect to a random measure.
\par
Every predictable and finite variation process $A$ decomposes as
\begin{equation}\label{eq:decomppredictcanon}
    A=A_{0}+A^{c}+x\ast\mu^{A},
\end{equation}
where $A^{c}$ is continuous and finite variation, and $A^{c}_{0}=0$ (see Proposition 3.15, \cite{pure}).
\section{Main results}\label{sec:main}
The main result of this section, Theorem \ref{thm:amaze}, is exceptionally surprising. It allows one to pass from a bounded sequence $\{X^{n}\}_{n}$ of semimartingales to a semimartingale limit in a $\sigma$-localized (in the sense of Kallsen \cite{sigmaloc}) version of the semimartingale topology.
\begin{theorem}\label{thm:amaze}
    Let $\{X^{n}\}_{n}$ be a sequence of semimartingales such that
    \begin{equation*}
        \mathrm{co}\left\{\vert{(\xi\cdot X^{n})_{T}\vert}:n\in\mathbb{N},\xi\in\mathscr{P}(1)\right\}
    \end{equation*}
    is bounded in probability.  Then there exists an increasing sequence $\{D_{n}\}_{n}$ of predictable sets increasing to $\Omega\times[0,T]$, a semimartingale $Y$, and $Y^{n}\in\mathrm{co}\{X^{m}:m\geq n\}$, such that $\{(\mathbf{1}_{D_{m}}\cdot Y^{n})\}_{n}$ converges to $(\mathbf{1}_{D_{m}}\cdot Y)$ in the semimartingale topology for each $m\in\mathbb{N}$.
\end{theorem}
The boundedness condition in the statement of Theorem \ref{thm:amaze} is a convex version of Stricker's notion of predictable uniform tightness (see \cite{embound}). Specific examples, and non-examples, of when the condition is satisfied are given in Section \ref{sec:super}. Here we highlight only that the condition is satisfied by any u.c.p. convergent sequence of supermartingales, at least after passing to a subsequence.
\par
We deduce Theorem \ref{thm:amaze} as a consequence of a slightly more general result, Theorem \ref{thm:amazedisjoint}, whose formulation `disjointifies' that of Theorem \ref{thm:amaze}. Its statement is similar to the $H^{1}$-Kadec-Pełczyński decomposition of Delbaen and Schachermayer (see Theorem C, \cite{del-sch}).
\begin{theorem}\label{thm:amazedisjoint}
    Let $\{X^{n}\}_{n}$ be a sequence of semimartingales such that
    \begin{equation*}
        \mathrm{co}\left\{\vert{(\xi\cdot X^{n})_{T}\vert}:n\in\mathbb{N},\xi\in\mathscr{P}(1)\right\}
    \end{equation*}
    is bounded in probability. Then there exists a disjoint sequence $\{G_{n}\}_{n}$ of predictable sets, a semimartingale $Y$, and $Y^{n}\in\mathrm{co}\{X^{m}:m\geq n\}$, such that $\{(\mathbf{1}_{(\Omega\times[0,T])\setminus G_{n}}\cdot Y^{n})\}_{n}$ converges to $Y$ in the semimartingale topology.
\end{theorem}
Theorem \ref{thm:amazedisjoint} implies Theorem \ref{thm:amaze}. Indeed, suppose $\{Y^{n}\}_{n}$, $Y$, and $\{G_{n}\}_{n}$ are obtained from $\{X^{n}\}_{n}$ using Theorem \ref{thm:amazedisjoint}. Defining $D_{n}=\bigcap_{m=n}^{\infty}(\Omega\times[0,T])\setminus G_{m}$ for each $n$, the sequence $\{D_{n}\}_{n}$ has the following properties:
    \begin{enumerate}
        \item $\{D_{n}\}_{n}$ is increasing, and $\Omega\times[0,T]=\bigcup_{n=1}^{\infty}D_{n}$.
        \item For each $m\in\mathbb{N}$, $\{(\mathbf{1}_{D_{m}}\cdot Y^{n})\}_{n}$ converges in the semimartingale topology to $(\mathbf{1}_{D_{m}}\cdot Y)$.
    \end{enumerate}
\section{Proof of Theorem \ref{thm:amazedisjoint}}\label{sec:proofs}
The proof of Theorem \ref{thm:amazedisjoint} progresses as follows. In \S\ref{subsec:decomp} and \S\ref{subsec:proofdecomp} we state and prove a decomposition result for sequences of semimartingales satisfying a boundedness condition (distinct from that of Theorem \ref{thm:amazedisjoint}). In \S\ref{subsec:lem} we return to the context of Theorem \ref{thm:amazedisjoint}, establishing a rather subtle relation between the boundedness conditions in \S\ref{subsec:decomp} and the assumptions of Theorem \ref{thm:amazedisjoint}. Finally, \S\ref{subsec:proving} finishes the proof.
\subsection{A decomposition result}\label{subsec:decomp}
The main technical difficulty in the proof of the fundamental theorem of asset pricing by Delbaen and Schachermayer \cite{ftap} is establishing an ad hoc convex compactness principle for a certain sequence $\{X^{n}\}_{n}$ in the semimartingale topology. Indeed, the majority of the proof is dedicated to finding a probability measure $\mathbb{Q}\sim\mathbb{P}$ such that each $X^{n}$ is a special semimartingale under $\mathbb{Q}$, and one can find $\widetilde{X}^{n}\in\mathrm{co}\{X^{m}:m\geq n\}$ such that the sequence of local martingale parts of the $\mathbb{Q}$-Doob-Meyer decompositions of the $\widetilde{X}^{n}$'s converge in the semimartingale topology.
\par
We explore a generalization of this idea. The main result of this subsection states that, if a sequence $\{X^{n}\}_{n}$ of semimartingales satisfies an appropriate boundedness condition, then we may find a probability measure $\mathbb{Q}\sim\mathbb{P}$ such that each $X^{n}$ is a special semimartingale under $\mathbb{Q}$, and the sequence of local martingale parts $\{M^{n}\}_{n}$ of the $\mathbb{Q}$-Doob-Meyer decomposition of $\{X^{n}\}_{n}$ converges in the semimartingale topology.
\begin{theorem}\label{thm:technicalmemin}
    Let $\{X^{n}\}_{n}$ be a sequence of semimartingales such that
    \begin{equation}\label{eq:l0quad}
        \mathrm{co}\left\{[X^{n},X^{n}]_{T}:n\in\mathbb{N}\right\},
    \end{equation}
    is bounded in probability. Then there exists a probability measure $\mathbb{Q}\sim\mathbb{P}$, and $\widetilde{X}^{n}\in\mathrm{co}\{X^{m}:m\geq n\}$ such that the following holds. $\widetilde{X}^{n}$ is a special semimartingale under $\mathbb{Q}$ for each $n$, and the local martingale parts $\{M^{n}\}_{n}$ of the $\mathbb{Q}$-Doob-Meyer decomposition $\widetilde{X}^{n}=\widetilde{X}^{n}_{0}+M^{n}+A^{n}$ of $\{\widetilde{X}^{n}\}_{n}$ converge in the semimartingale topology.
    \par
    If, furthermore,
    \begin{equation}\label{eq:lastcond}
        \mathrm{co}\left\{\vert{(\xi\cdot X^{n})_{T}\vert}:n\in\mathbb{N},\xi\in\mathscr{P}(1)\right\},
    \end{equation}
    is bounded in probability, then $\mathbb{Q}$ as above may be chosen so that
    \begin{equation*}
        \sup_{n}\int_{\Omega}\mathrm{var}(A^{n})_{T}d\mathbb{Q}<\infty.
    \end{equation*}
\end{theorem}
Similar results have been obtained by Mémin \cite{memin}, albeit under the assumption that the sequence $\{X^{n}\}_{n}$ converges in the semimartingale topology—a severe limitation in our context.
\subsection{Proof of Theorem \ref{thm:technicalmemin}}\label{subsec:proofdecomp}
In the proof of Theorem \ref{thm:technicalmemin}, it is necessary to switch the probability measure. As such, for much of this subsection, we work with an equivalent probability measure $\mathbb{Q}\sim\mathbb{P}$.
\begin{lemma}\label{lem:weakcompactness}
    Let $\{M^{n}\}_{n}$ be a sequence of $\mathbb{Q}$-local martingales starting at $M^{n}_{0}=0$. If
    \begin{equation*}
        \sup_{n}\int_{\Omega}\left[M^{n},M^{n}\right]_{T}d\mathbb{Q}<\infty,
    \end{equation*}
    then there exists $N^{n}\in\mathrm{co}\{X^{m}:m\geq n\}$ such that $\{N^{n}\}_{n}$ converges in the semimartingale topology. Furthermore, if $L^{n}\in\mathrm{co}\{N^{m}:m\geq n\}$, then $\{L^{n}\}_{n}$ also converges in the semimartingale topology to the same limit.
\end{lemma}
\begin{proof}
    Clearly, each $M^{n}$ is a martingale (indeed, each $M^{n}$ has integrable quadratic variation). Let $\mathfrak{H}$ denote the vector space of martingales $M$ starting at zero such that
    \begin{equation*}
        \int_{\Omega}[M,M]_{T}d\mathbb{Q}<\infty.
    \end{equation*}
    Then $\mathfrak{H}$ is a Hilbert space under the inner product $(M,N)\longmapsto\int_{\Omega}[M,N]_{T}d\mathbb{Q}$. By the Bourbaki-Alaouglu theorem, any bounded subset of $\mathfrak{H}$ is relatively weakly compact. Thus, the Eberlein-Šmulian theorem (see Theorem 3.19, \cite{brezis}) and Mazur's lemma (see Corollary 3.8, \cite{brezis}) implies the existence of $N^{n}\in\mathrm{co}\left\{M^{m}:m\geq n\right\}$ such that $\{N^{n}\}_{n}$ converges in the norm of $\mathfrak{H}$; it follows that $\{N^{n}\}_{n}$ converges in the semimartingale topology. Indeed, we may assume (via translation) that $\{N^{n}\}_{n}$ converges to zero in norm, implying that for arbitrary $\varepsilon>0$,
    \begin{equation*}
        \sup_{\xi\in\mathscr{P}(1)}\mathbb{Q}\left(\left\{\left(\xi\cdot N^{n}\right)^{\ast}_{T}\geq\varepsilon\right\}\right)\leq\frac{1}{\varepsilon^{2}}\sup_{\xi\in\mathscr{P}(1)}\int_{\Omega}\left(\left(\xi\cdot N^{n}\right)^{\ast}_{T}\right)^{2}d\mathbb{Q}
    \end{equation*}
    \begin{equation}\label{eq:quadleq}
        \leq \frac{C}{\varepsilon^{2}}\sup_{\xi\in\mathscr{P}(1)}\int_{\Omega}\int_{0}^{T}\xi^{2}d[N^{n},N^{n}]_{T}d\mathbb{Q}\leq\frac{C}{\varepsilon^{2}}\int_{\Omega}[N^{n},N^{n}]_{T}d\mathbb{Q},
    \end{equation}
    for some universal constant $C>0$ by Markov's concentration inequality and the Burkholder-Davis-Gundy inequality. Taking $n$ to infinity in the right-hand side of (\ref{eq:quadleq}) yields the claim. 
    \par
    For the last part, simply note that convergence in the norm of $\mathfrak{H}$ is stable under passing to forward convex combinations.
\end{proof}
\begin{remark}
    More generally than above, Lemma \ref{lem:weakcompactness} holds whenever the set of square functions $\left\{[M^{n},M^{n}]_{T}^{\frac{1}{2}}:n\in\mathbb{N}\right\}$ is uniformly $\mathbb{Q}$-integrable, essentially as a consequence of the Dellacherie-Meyer-Yor characterization of weak compactness in the martingale Hardy space $H^{1}(\mathbb{Q})$ (see \cite{delmeyyor}). We have singled out the $L^{2}$-bounded case since it is, in some sense, the most fundamental case (c.f. Théorème II.3, \cite{memin}).
\end{remark}
We will need the following lemma about the convex hull of finite unions of convex bounded subsets of the nonnegative cone in $L^{0}$.
\begin{lemma}\label{lem:union}
    Let $K$ and $L$ be convex sets of nonnegative random variables which are bounded in probability. Then $\mathrm{co}(K\cup L)$ is bounded in probability.
\end{lemma}
\begin{proof}
    By (Lemma 2.3, \cite{bipolar}), we may find equivalent probability measures $\mathbb{Q}_{1}\sim\mathbb{P}$, $\mathbb{Q}_{2}\sim\mathbb{P}$, such that
    \begin{equation*}
        \sup_{\xi\in K}\int_{\Omega}\xi d\mathbb{Q}_{1}<\infty,
    \end{equation*}
    \begin{equation*}
        \sup_{\xi\in L}\int_{\Omega}\xi d\mathbb{Q}_{2}<\infty.
    \end{equation*}
    Define a new probability measure $\mathbb{Q}\sim\mathbb{P}$ by its Radon-Nikodým derivative with respect to $\mathbb{P}$,
    \begin{equation*}
        \frac{d\mathbb{Q}}{d\mathbb{P}}=\frac{1}{C}\left(\frac{d\mathbb{Q}_{1}}{d\mathbb{P}}\wedge\frac{d\mathbb{Q}_{2}}{d\mathbb{P}}\right),
    \end{equation*}
    where $C=\int_{\Omega}\left(\frac{d\mathbb{Q}_{1}}{d\mathbb{P}}\wedge\frac{d\mathbb{Q}_{2}}{d\mathbb{P}}\right)d\mathbb{P}$ clearly satisfies $0<C<\infty$. Then
    \begin{equation*}
        \sup_{\xi\in K\cup L}\int_{\Omega}\xi d\mathbb{Q}\leq\frac{1}{C}\left(\left(\sup_{\xi\in K}\int_{\Omega}\xi d\mathbb{Q}_{1}\right)\vee\left(\sup_{\xi\in L}\int_{\Omega}\xi d\mathbb{Q}_{2}\right)\right)<\infty,
    \end{equation*}
    and this inequality extends to $\mathrm{co}(K\cup L)$ by the triangle inequality. Thus, $\mathrm{co}(K\cup L)$ is bounded in $L^{1}(\mathbb{Q})$ for some probability measure $\mathbb{Q}\sim\mathbb{P}$, which implies boundedness in probability.
\end{proof}
We are now ready to prove Theorem \ref{thm:technicalmemin}.
\begin{proof}
Define
\begin{equation*}
    K=\mathrm{co}\left\{[X^{n},X^{n}]_{T}:n\in\mathbb{N}\right\}.
\end{equation*}
Our assumptions imply that $K$ is bounded in probability. By (Lemma 2.3, \cite{bipolar}), we therefore may find a probability measure $\mathbb{Q}\sim\mathbb{P}$ such that
    \begin{equation}\label{eq:l1bound}
        \sup_{\xi\in K}\int_{\Omega}\xi d\mathbb{Q}<\infty.
    \end{equation}
Clearly, $X^{n}$ is a special semimartingale under $\mathbb{Q}$ for each $n$, as $[X^{n},X^{n}]_{T}\in L^{1}(\mathbb{Q})$ for each $n$. Denote by $X^{n}=X^{n}_{0}+N^{n}+B^{n}$ its Doob-Meyer decomposition with respect to $\mathbb{Q}$, where $B^{n}$ is the compensator process. We have the inequality
\begin{equation*}
    \int_{\Omega}[X^{n},X^{n}]_{T}d\mathbb{Q}\geq\int_{\Omega}[N^{n},N^{n}]_{T}d\mathbb{Q},
\end{equation*}
for each $n$, and so $\sup_{n}\int_{\Omega}[N^{n},N^{n}]_{T}d\mathbb{Q}<\infty$ (see (\ref{eq:l1bound})). Consequently, we may apply Lemma \ref{lem:weakcompactness} to obtain $M^{n}\in\mathrm{co}\{N^{m}:m\geq n\}$ such that $\{M^{n}\}_{n}$ converges in the semimartingale topology. Denote by $\widetilde{X}^{n}$ the corresponding convex combination of $\{X_{n},X_{n+1},\dots\}$, and by $A^{n}$ the corresponding convex combination of $\{B^{n},B^{n+1},\dots\}$. Since $\widetilde{X}^{n}=\widetilde{X}^{n}_{0}+M^{n}+A^{n}$ is the Doob-Meyer decomposition of $\widetilde{X}^{n}$, and $\{M^{n}\}_{n}$ converges in the semimartingale topology, this yields the first part of the claim.
\par
Suppose now that
\begin{equation*}
    J=\mathrm{co}\left\{\vert{(\xi\cdot X^{n})_{T}\vert}:n\in\mathbb{N},\xi\in\mathscr{P}(1)\right\}
\end{equation*}
is bounded in probability. In light of Lemma \ref{lem:union} and (Lemma 2.3, \cite{bipolar}), we may assume that $J$ is bounded in $L^{1}(\mathbb{Q})$. It suffices to show that $\{\mathrm{var}(A^{n})\}_{n}$ is bounded in $L^{1}(\mathbb{Q})$. For each $n$, there exists a predictable process $\xi^{n}\in\mathscr{P}(1)$ with $\mathrm{var}(A^{n})=(\xi^{n}\cdot A^{n})$. Thus, by the triangle inequality,
    \begin{equation*}
        \mathrm{var}(A^{n})_{T}\leq\vert{(\xi^{n}\cdot \widetilde{X}^{n})_{T}\vert}+\vert{(\xi^{n}\cdot M^{n})_{T}\vert},
    \end{equation*}
    from which we can conclude that
    \begin{equation*}
        \int_{\Omega}\mathrm{var}(A^{n})_{T}d\mathbb{Q}\leq\int_{\Omega}\vert{(\xi^{n}\cdot \widetilde{X}^{n})}\vert_{T}d\mathbb{Q}+\left(\int_{\Omega}(\xi^{n}\cdot M^{n})_{T}^{2}d\mathbb{Q}\right)^{\frac{1}{2}}
    \end{equation*}
    \begin{equation*}
        \leq\int_{\Omega}\vert{(\xi^{n}\cdot \widetilde{X}^{n})}\vert_{T}d\mathbb{Q}+\left(\int_{\Omega}[M^{n},M^{n}]_{T}d\mathbb{Q}\right)^{\frac{1}{2}}\leq C<\infty,
    \end{equation*}
    for some absolute constant $C>0$ by Itô's isometry.
\end{proof}
\subsection{Some lemmata}\label{subsec:lem}
We now return to the context of Theorem \ref{thm:amazedisjoint}. Recall that $\{X^{n}\}_{n}$ denotes a sequence of semimartingales such that
\begin{equation*}
        \mathrm{co}\left\{\vert{(\xi\cdot X^{n})_{T}\vert}:n\in\mathbb{N},\xi\in\mathscr{P}(1)\right\}
\end{equation*}
is bounded in probability. Since both the assumption and conclusion of Theorem \ref{thm:amazedisjoint} are valid for $\{X^{n}\}_{n}$ iff they are valid for $\{X^{n}-X^{n}_{0}\}$, without loss of generality we will assume up until the end of the proof of Theorem \ref{thm:amazedisjoint} that $X^{n}_{0}=0$ for all $n\in\mathbb{N}$.
\par
We begin with the following observation.
\begin{lemma}\label{lem:tyyor}
    The set
    \begin{equation}\label{eq:l0max}
        \mathrm{co}\left\{(X^{n})^{\ast}_{T}:n\in\mathbb{N}\right\},
    \end{equation}
    is bounded in probability.
\end{lemma}
\begin{proof}
    By (Lemma 2.3, \cite{bipolar}), there exists an equivalent probability measure $\mathbb{Q}\sim\mathbb{P}$ such that
    \begin{equation*}
        \sup_{n}\sup_{\xi\in\mathscr{P}(1)}\int_{\Omega}\vert{(\xi\cdot X^{n})_{T}}\vert d\mathbb{Q}<\infty.
    \end{equation*}
    By a result of Yor (see e.g. Theorem 104 in Chapter VII of \cite{dellmey}) combined with equation (98.7) from Chapter VII of \cite{dellmey} yields that
    \begin{equation*}
        \sup_{n}\int_{\Omega}\left(X^{n}\right)^{\ast}_{T}d\mathbb{Q}<\infty,
    \end{equation*}
    and so the claim follows from Markov's inequality and convexity.
\end{proof}
We need the following lemma to apply the results of \S\ref{subsec:decomp}. Similar conclusions have been obtained under predictable uniform tightness, but without the assumption and conclusion of convex boundedness (see Lemme 1.2, \cite{ptrfskoro}).
\begin{lemma}\label{lem:removerestrict}
    We may find convex combinations $Y^{n}\in\mathrm{co}\{X^{m}:m\geq n\}$
    \begin{equation*}
        \mathrm{co}\left\{[Y^{n},Y^{n}]_{T}:n\in\mathbb{N}\right\},
    \end{equation*}
    is bounded in probability.
\end{lemma}
Remark that the unintuitive passage from $\{X^{n}\}_{n}$ to convex combinations $\{Y^{n}\}_{n}$ is actually necessary in the context of Lemma \ref{lem:removerestrict}. We now give a sketch of a counterexample, showing that in general $\mathrm{co}\left(\left\{[X^{n},X^{n}]_{T}:n\in\mathbb{N}\right\}\right)$ may be unbounded under the assumptions of Theorem \ref{thm:amaze}.
\begin{example}
    Suppose that $T=1$, and $(\Omega,\mathscr{F}_{1/2},\mathbb{P})$ is a Lebesgue-Rokhlin probability space. Gao, Leung, and Xanthos \cite{gaoexample} provide a construction which in the present context implies the existence of an $L^{0}$-null $L^{1}(\mathbb{P})$-bounded sequence $\{g_{n}\}_{n}$ of $\mathscr{F}_{1/2}$-measurable random variables such that $\{g_{n}\}_{n}$ is not uniformly $\mathbb{Q}$-integrable under any probability measure $\mathbb{Q}\sim\mathbb{P}$. Define $\{X^{n}\}_{n}$ by
\begin{equation*}
    X^{n}=\mathbf{1}_{\llbracket{1/2,1}\rrbracket}g_{n}.
\end{equation*}
Evidently
\begin{equation*}
    [X^{n},X^{n}]_{1}=g_{n}^{2},
\end{equation*}
and $\{X^{n}\}_{n}$ satisfies the hypotheses of Theorem \ref{thm:amazedisjoint}. Suppose that $\mathrm{co}\{[X^{n},X^{n}]_{1}:n\in\mathbb{N}\}$ is bounded in probability; by the above identity in combination with (Lemma 2.3, \cite{bipolar}), we obtain the existence of an equivalent probability measure $\mathbb{Q}\sim\mathbb{P}$ such that $\{g_{n}\}_{n}$ is bounded in $L^{2}(\mathbb{Q})$. Thus, $\{g_{n}\}_{n}$ is uniformly $\mathbb{Q}$-integrable, a contradiction. 
\end{example}
\par
To prove Lemma \ref{lem:removerestrict}, we use the following elementary inequality.
\begin{lemma}\label{lem:easy}
    Let $\mathbb{P}$ be any probability measure. Suppose that a random variable $\xi\geq0$ is such that
    \begin{equation*}
        \xi\leq\sum_{i=1}^{n}\xi_{i},
    \end{equation*}
    where $\xi_{i}\geq0$ are random variables. Then
    \begin{equation*}
        \mathbb{P}(\{\xi>\eta\})\leq\sum_{i=1}^{n}\mathbb{P}\left(\left\{\xi_{i}>\frac{\eta}{n}\right\}\right),
    \end{equation*}
    for any $\eta\geq0$.
\end{lemma}
\begin{proof}
    Clearly, $\{\xi>\eta\}\subset\bigcup_{i=1}^{n}\left\{\xi_{i}>\frac{\eta}{n}\right\}$. Thus,
    \begin{equation*}
        \mathbb{P}\left(\{\xi>\eta\}\right)\leq\mathbb{P}\left(\bigcup_{i=1}^{n}\left\{\xi_{i}>\frac{\eta}{n}\right\}\right)\leq\sum_{i=1}^{n}\mathbb{P}\left(\left\{\xi_{i}>\frac{\eta}{n}\right\}\right),
    \end{equation*}
    by the union bound, as desired.
\end{proof}
We will shortly begin the proof of Lemma \ref{lem:removerestrict}. To this end, denote by $\Delta_{\infty}$ the following infinite-dimensional simplex:
\begin{equation*}
    \Delta_{\infty}=\left\{a\in[0,1]^{\mathbb{N}}:\sum_{n=1}^{\infty}a_{n}=1,\textrm{ and there exists }m\in\mathbb{N}\textrm{ with }a_{n}=0\textrm{ for all $n\geq m$}\right\}.
\end{equation*}
\begin{proof}[Proof of Lemma \ref{lem:removerestrict}]
In light of Lemma \ref{lem:tyyor}, the set
    \begin{equation*}
        \mathrm{co}\left\{\left(X^{n}\right)^{\ast}_{T}:n\in\mathbb{N}\right\}
    \end{equation*}
    is bounded in probability. By (Lemma 2.3, \cite{bipolar}) and Lemma \ref{lem:union}, we may find an equivalent probability measure $\mathbb{Q}\sim\mathbb{P}$ such that
    \begin{equation*}
        \sup_{\xi\in\mathscr{P}(1)}\sup_{n}\int_{\Omega}\vert{(\xi\cdot X^{n})_{T}}\vert d\mathbb{Q}=C_{1}<\infty,
    \end{equation*}
    \begin{equation*}
        \sup_{n}\int_{\Omega}(X^{n})^{\ast}_{T}d\mathbb{Q}=C_{2}<\infty.
    \end{equation*}
    By (Lemma 2.5, \cite{schruessdies}), we may find a random variable $\zeta>0$, and convex combinations $Y^{n}\in\mathrm{co}\{X^{m}:m\geq n\}$, such that
    \begin{equation*}
        (Y^{n})^{\ast}_{T}\leq\zeta,
    \end{equation*}
    for each $n$.
    \par
    Let $\varepsilon>0$. It suffices to show that there exists $K>0$ such that for any $\lambda\in\Delta_{\infty}$, we have that
    \begin{equation*}
        \mathbb{Q}\left(\left\{\sum_{i}\lambda_{i}[Y^{i},Y^{i}]_{T}>K\right\}\right)<\varepsilon.
    \end{equation*}
    We may find $K_{1}\geq1$ such that $\mathbb{Q}\left(\left\{\zeta^{2}>K_{1}\right\}\right)<\frac{\varepsilon}{3}$. By Markov's inequality, we have that
    \begin{equation}\label{eq:markovint}
        \sup_{\{\xi^{i}\}_{i}\subset\mathscr{P}(1)}\sup_{n}\mathbb{Q}\left(\left\{\sum_{i}\lambda_{i}\vert{(\xi^{i}\cdot Y^{i})_{T}\vert}>\eta\right\}\right)\leq\frac{C_{1}}{\eta},
    \end{equation}
    for any $\eta>0$ and any $\lambda\in\Delta_{\infty}$. Let $K_{2}>0$ be such that $\frac{2C_{1}K_{1}}{K_{2}}<\frac{\varepsilon}{3}$.
    \par
    Define $K=3(K_{1}\vee K_{2})$ and $\widetilde{K}=K_{1}\vee K_{2}=\frac{K}{3}$. For each $\lambda\in\Delta_{\infty}$ and $\eta\geq0$, we may (using integration by parts) write
    \begin{equation*}
        \sum_{i}\lambda_{i}[Y^{i},Y^{i}]_{T}=\sum_{i}\lambda_{i}\left((Y^{i}_{T})^{2}-2(\mathbf{1}_{\{\vert{Y^{i}_{-}\vert}\leq\eta\}}Y^{i}_{-}\cdot Y^{i})_{T}-2(\mathbf{1}_{\{\vert{Y^{i}_{-}\vert}>\eta\}}Y^{i}_{-}\cdot Y^{i})_{T}\right).
    \end{equation*}
    Thus, by Lemma \ref{lem:easy} and (\ref{eq:markovint}),
    \begin{equation*}
        \mathbb{Q}\left(\left\{\sum_{i}\lambda_{i}[Y^{i},Y^{i}]_{T}>K\right\}\right)\leq\mathbb{Q}\left(\left\{\zeta^{2}>\widetilde{K}\right\}\right)+\mathbb{Q}\left(\left\{\sum_{i}\lambda_{i}2\vert(\mathbf{1}_{\{\vert{Y^{i}_{-}\vert}\leq K_{1}\}}Y^{i}_{-}\cdot Y^{i})_{T}\vert>\widetilde{K}\right\}\right)
    \end{equation*}
    \begin{equation*}
        +\mathbb{Q}\left(\left\{\sum_{i}\lambda_{i}2\vert(\mathbf{1}_{\{\vert{Y^{i}_{-}\vert}>K_{1}\}}Y^{i}_{-}\cdot Y^{i})_{T}\vert>\widetilde{K}\right\}\right)\leq2\mathbb{Q}\left(\left\{\zeta^{2}>K_{1}\right\}\right)
    \end{equation*}
    \begin{equation*}
        +\mathbb{Q}\left(\left\{\sum_{i}\lambda_{i}\vert(\frac{1}{K_{1}}\mathbf{1}_{\{\vert{Y^{i}_{-}\vert}\leq K_{1}\}}Y^{i}_{-}\cdot Y^{i})_{T}\vert>\frac{K_{2}}{2K_{1}}\right\}\right)<\frac{2\varepsilon}{3}+\frac{2C_{1}K_{1}}{K_{2}}<\varepsilon,
    \end{equation*}
    as desired.
\end{proof}
\subsection{The proof}\label{subsec:proving}
In this subsection, we complete the proof of Theorem \ref{thm:amazedisjoint}. The final technical lemma needed is the following version of the Kadec-Pełczyński decomposition, first established by \cite{ogkp}.
\begin{lemma}\label{lem:ogkp}
    Let $(E,\Sigma,\mu)$ be a finite measure space (not necessarily $\mu$-complete). If $\{\xi_{n}\}_{n}\subset L^{1}(\mu)$ satisfies
    \begin{equation*}
        \sup_{n}\int_{E}\vert{\xi_{n}}\vert d\mu<\infty,
    \end{equation*}
    then there exists a subsequence $\{n_{k}\}_{k}$ and a disjoint sequence $\{G_{k}\}_{k}\subset\Sigma$ such that $\{\mathbf{1}_{X\setminus G_{k}}\xi_{n_{k}}\}_{k}$ is uniformly $\mu$-integrable.
\end{lemma}
\begin{proof}
    See (Lemma 5.2.8, \cite{albkalton}).
\end{proof}
We are now ready to give the proof of Theorem \ref{thm:amazedisjoint}. First, we sketch out a rough outline. Modulo some technicalities, Theorem \ref{thm:technicalmemin} allows one to switch to $\mathbb{Q}\sim\mathbb{P}$ under which all the $X^{n}$'s are special, and under which the martingale parts $\{M^{n}\}_{n}$ of the $\mathbb{Q}$-Doob-Meyer decompositions $X^{n}=M^{n}+A^{n}$ converge in the semimartingale topology. Thus, we may focus entirely on the sequence $\{A^{n}\}_{n}$ of compensators, which we deal with via the Kadec-Pełczyński theorem (see above) applied to an $L^{1}$-space on the predictable $\sigma$-algebra. This $L^{1}$-space is constructed using a superposition of the sequence of Doléans measures obtainable from $\{\mathrm{var}(A^{n})\}_{n}$.
\begin{proof}[Proof of Theorem \ref{thm:amaze}]
    Applying Lemma \ref{lem:removerestrict} and Theorem \ref{thm:technicalmemin}, we may find a probability measure $\mathbb{Q}\sim\mathbb{P}$ and convex combinations $Y^{n}\in\mathrm{co}\{X^{m}:m\geq n\}$ such that each $Y^{n}$ is a special semimartingale under $\mathbb{Q}$, $Y^{n}=M^{n}+A^{n}$ is the $\mathbb{Q}$-Doob-Meyer decomposition of $Y^{n}$, $\{M^{n}\}_{n}$ converges in the semimartingale topology to some $M$, and
    \begin{equation}\label{eq:predictablel1bound}
        C=\sup_{n}\int_{\Omega}\mathrm{var}(A^{n})_{T}d\mathbb{Q}<\infty.
    \end{equation}
    Our goal will be to apply the Kadec-Pełczyński-type decomposition provided by Lemma \ref{lem:ogkp} to a measure space constructed on $(\Omega\times[0,T],\mathscr{P})$, where $\mathscr{P}$ denotes the predictable $\sigma$-algebra. Define a measure $\nu:\mathscr{P}\longrightarrow[0,\infty]$ by
    \begin{equation*}
        \nu=\left(D\longmapsto\sum_{n=1}^{\infty}\frac{1}{2^{n}}\int_{\Omega}(\mathbf{1}_{D}\cdot\mathrm{var}(A^{n}))_{T}d\mathbb{Q}\right).
    \end{equation*}
    It is easy to see that $\nu$ is a finite measure, as $\nu(\Omega\times[0,T])\leq C$.
    \par
    There is a càdlàg, predictable, increasing, and integrable process $V$ with
    \begin{equation*}
        \nu(D)=\int_{\Omega}(\mathbf{1}_{D}\cdot V)_{T}d\mathbb{Q},
    \end{equation*}
    for any $D\in\mathscr{P}$ (see p. 128, \cite{dellmey}). We may use the predictable Radon-Nikodým theorem established by Delbaen and Schachermayer (see Theorem 2.1(ii), \cite{predictablernderivative}) to find, for each $n$, a predictable process $\varphi^{n}\in L^{1}(\nu)$ and a predictable set $C^{n}$ with
    \begin{equation*}
        A^{n}=(\varphi^{n}\cdot V)+(\mathbf{1}_{C^{n}}\cdot A^{n}),
    \end{equation*}
    \begin{equation*}
        (\mathbf{1}_{C^{n}}\cdot V)=0.
    \end{equation*}
    However, the condition that $(\mathbf{1}_{C^{n}}\cdot V)=0$ clearly also implies that $(\mathbf{1}_{C^{n}}\cdot A^{n})=0$; indeed,
    \begin{equation*}
        \int_{\Omega}\vert{(\mathbf{1}_{C^{n}}\cdot A^{n})_{t}}\vert d\mathbb{Q}\leq\int_{\Omega}(\mathbf{1}_{C^{n}}\cdot\mathrm{var}(A^{n}))_{t}d\mathbb{Q}
    \end{equation*}
    \begin{equation*}
        \leq 2^{n}\nu(C^{n})=2^{n}\int_{\Omega}(\mathbf{1}_{C^{n}}\cdot V)_{T}d\mathbb{Q}=0,
    \end{equation*}
    so that $(\mathbf{1}_{C^{n}}\cdot A^{n})_{t}=0$ almost surely for each $t\in[0,T]$, so we may apply (Problem 1.5, \cite{ks91}) to conclude that $(\mathbf{1}_{C^{n}}\cdot A^{n})=0$ up to evanescence. Thus, $A^{n}=(\varphi^{n}\cdot V)$ and $\mathrm{var}(A^{n})=(\vert{\varphi^{n}}\vert\cdot V)$.
    \par
    It is easy to see that
    \begin{equation*}
        \sup_{n}\int_{\Omega\times[0,T]}\vert{\varphi^{n}}\vert d\nu=C<\infty,
    \end{equation*}
    where $C>0$ is the constant from (\ref{eq:predictablel1bound}). By Komlós's theorem (see Theorem 1.3, \cite{del-sch}), we therefore may pass to convex combinations (still denoted $\{\varphi^{n}\}_{n}$) converging in the following sense. There exists a predictable set $H$ such that
    \begin{equation*}
        \lim_{n}\varphi_{n}=\varphi,
    \end{equation*}
    pointwise on $H$ for some predictable process $\varphi\in L^{1}(\nu)$, and $\nu((\Omega\times[0,T])\setminus H)=0$. Our unchanged notation is justified by the preservation under convex combinations of all properties relevant to the sequel (including semimartingale convergence of $\{M^{n}\}_{n}$ to $M$, in light of Lemma \ref{lem:weakcompactness} and the proof of Theorem \ref{thm:technicalmemin}). By Lemma \ref{lem:ogkp} and passing to a subsequence if necessary, there exists a disjoint sequence $\{G_{n}\}_{n}$ of predictable sets such that $\{\mathbf{1}_{(\Omega\times[0,T])\setminus G_{n}}\varphi^{n}\}_{n}$ is uniformly integrable (with respect to $\nu$). By Vitali's convergence theorem, we have that,
    \begin{equation*}
        \lim_{n}\int_{\Omega\times[0,T]}\vert{\mathbf{1}_{(\Omega\times[0,T])\setminus G_{n}}\varphi^{n}-\varphi}\vert d\nu=0.
    \end{equation*}
    Thus,
    \begin{equation*}
        \lim_{n}\int_{\Omega}\mathrm{var}((\mathbf{1}_{(\Omega\times[0,T])\setminus G_{n}}\varphi^{n}\cdot V)-(\varphi\cdot V))_{T}d\mathbb{Q}=0,
    \end{equation*}
    and hence Markov's inequality implies $\{\mathrm{var}((\mathbf{1}_{(\Omega\times[0,T])\setminus G_{n}}\varphi^{n}\cdot V)-(\varphi\cdot V))_{T}\}_{n}$ converges to zero in probability. In particular, $\{(\mathbf{1}_{(\Omega\times[0,T])\setminus G_{n}}\varphi^{n}\cdot V)\}_{n}$ converges in the semimartingale topology to $(\varphi \cdot V)$ (see Proposition 2.7, \cite{kardem}). Since we may write $(\mathbf{1}_{(\Omega\times[0,T])\setminus G_{n}}\varphi^{n}\cdot V)=(\mathbf{1}_{(\Omega\times[0,T])\setminus G_{n}}\cdot A^{n})$, it follows that $\{(\mathbf{1}_{(\Omega\times[0,T])\setminus G_{n}}\cdot A^{n})\}_{n}$ converges in the semimartingale topology to $(\varphi\cdot V)$.
    \par
    Since $\{M^{n}\}_{n}$ converges in the semimartingale topology to $M$, $\{(\mathbf{1}_{(\Omega\times[0,T])\setminus G_{n}}\cdot M^{n})\}_{n}$ converges in the semimartingale topology to $M$. Thus, $(\mathbf{1}_{(\Omega\times[0,T])\setminus G_{n}}\cdot Y^{n})=(\mathbf{1}_{(\Omega\times[0,T])\setminus G_{n}}\cdot A^{n})+(\mathbf{1}_{(\Omega\times[0,T])\setminus G_{n}}\cdot M^{n})$ is a sum of convergent sequences in the semimartingale topology; it follows that $\{(\mathbf{1}_{(\Omega\times[0,T])\setminus G_{n}}\cdot Y^{n})\}_{n}$ converges in the semimartingale topology to $Y$, where $Y=(\varphi\cdot V)+M$. This proves the claim, as $\{G_{n}\}_{n}$ is a disjoint sequence.
\end{proof}
\begin{remark}\label{rem:singularpart}
    The above proof implies that the `remainder' $\{(\mathbf{1}_{G_{n}}\cdot Y^{n})\}_{n}$ can be decomposed as $(\mathbf{1}_{G_{n}}\cdot Y^{n})=N^{n}+B^{n}$, where $\{N^{n}\}_{n}$ is a null sequence in the semimartingale topology, $\{B^{n}\}_{n}$ is predictable and finite variation, and $\sup_{n}\int_{\Omega}\mathrm{var}(B^{n})_{T}d\mathbb{Q}<\infty$ for some $\mathbb{Q}\sim\mathbb{P}$. Many convex compactness theorems apply to this situation (see e.g. Proposition 13 of \cite{schcamp}).
\end{remark}
\section{Applications to semimartingale convergence}\label{sec:emeryconvg}
The Dunford-Pettis theorem asserts the equivalence between the following for a sequence $\{f_{n}\}_{n}$ in $L^{1}(\mu)$ (see Theorem 5.2.9, \cite{albkalton}).
\begin{enumerate}
    \item There exists $g_{n}\in\mathrm{co}\{f_{m}:m\geq n\}$ such that $\{g_{n}\}_{n}$ converges in $L^{1}(\mu)$.
    \item There exists $h_{n}\in\mathrm{co}\{f_{m}:m\geq n\}$ such that
    \begin{equation*}
        \lim_{n}\sup_{m}\int_{G_{n}}\vert{h_{m}}\vert d\mu=0
    \end{equation*}
    for all disjoint sequences $\{G_{n}\}_{n}$ of measurable sets, and $\{h_{n}\}_{n}$ is bounded in $L^{1}(\mu)$.
\end{enumerate}
In this section, we establish a Dunford-Pettis-type theorem for semimartingales. It characterizes when one can pass from a sequence $\{X^{n}\}_{n}$ of semimartingales to convex combinations converging in the semimartingale topology. Such questions have received an ad hoc treatment in arbitrage theory due to their immense applicability (see e.g. \cite{ftap} or \cite{emeryconvcomp}). However, outside of the narrow confines of arbitrage theory, the conditions in the literature are rather artificial, and never necessary.
\begin{theorem}\label{thm:dunfpett}
    Let $\{X^{n}\}_{n}$ be a sequence of semimartingales. There exists $Y^{n}\in\mathrm{co}\{X^{m}:m\geq n\}$ such that $\{Y^{n}\}_{n}$ converges in the Émery topology if, and only if, the following holds. There exists $Z^{n}\in\mathrm{co}\{X^{m}:m\geq n\}$ such that:
    \begin{enumerate}
        \item The set
        \begin{equation*}
            \mathrm{co}\left\{\vert{Z^{n}_{0}}\vert+\vert{(\xi\cdot Z^{n})_{T}}\vert:n\in\mathbb{N},\xi\in\mathscr{P}(1)\right\},
        \end{equation*}
        is bounded in probability.
        \item For every disjoint sequence $\{G_{n}\}_{n}$ of predictable sets,
        \begin{equation*}
            \lim_{n}\sup_{Z\in\mathrm{co}\{Z^{m}\}_{m}}\mathbf{D}((\mathbf{1}_{G_{n}}\cdot Z),0)=0.
        \end{equation*}
    \end{enumerate}
\end{theorem}
Before starting the proof of Theorem \ref{thm:dunfpett}, let us make a technical note. The necessity of $(1)\wedge(2)$ for \textit{some} convex combination has a non-trivial proof, and in general an Émery-convergent convex combination $\{Y^{n}\}_{n}$ of a sequence $\{X^{n}\}_{n}$ can fail both $(1)$ and $(2)$.
\begin{example}
    Suppose $T=1$ and $(\Omega,\mathscr{F}_{1/2},\mathbb{P})$ is a Lebesgue-Rokhlin probability space. By (Example 1.2, \cite{kardzit}), there exists an $L^{0}$-null sequence $\{g_{n}\}_{n}$ of nonnegative $\mathscr{F}_{1/2}$-measurable random variables such that the $L^{0}$-closure of $\mathrm{co}\{g_{n}\}_{n\geq m}$ contains all nonnegative random variables for all $m$. Let $\{t_{n}\}_{n}\subset(1/2,1)$ be a strictly increasing sequence. Define, for each $n$, $X^{n}=\mathbf{1}_{\llbracket{t_{n},1}\rrbracket}g_{n}$; then $\{X^{n}\}_{n}$ satisfies neither (1) nor (2), but converges to zero in the semimartingale topology. Indeed, $\neg$(1) is obvious, while $\neg(2)$ can be proven as follows. There exists a disjoint sequence $\{F_{n}\}_{n}$ of finite subsets of $\mathbb{N}$ such that the following holds. There exists $\{\lambda^{n}\}_{n}\subset\Delta_{\infty}$ satisfying $\{\lambda^{n}\neq0\}\subset F_{n}$ such that
    \begin{equation*}
        \int_{\Omega}\left\vert{h_{n}-1}\right\vert\wedge1 d\mathbb{P}<\frac{1}{n},
    \end{equation*}
    where $h_{n}=\sum_{i}\lambda^{n}_{i}g_{i}$. Define a disjoint sequence $\{G_{n}\}_{n}\subset\mathscr{P}$ by $G_{n}=\bigcup_{m\in F_{n}}\llbracket{t_{m}}\rrbracket$. Then
    \begin{equation*}
        \sup_{Z\in\mathrm{co}\{X^{m}\}_{m}}\mathbf{D}((\mathbf{1}_{G_{n}}\cdot Z),0)\geq\sup_{Z\in\mathrm{co}\{X^{m}\}_{m}}\int_{\Omega}\vert{(\mathbf{1}_{G_{n}}\cdot Z)_{T}}\vert\wedge1 d\mathbb{P}
    \end{equation*}
    \begin{equation*}
        \geq\int_{\Omega}\left\vert{h_{n}}\right\vert\wedge1 d\mathbb{P}\geq1-\frac{1}{n}\geq\frac{1}{2},
    \end{equation*}
    for $n\geq 2$. In particular, (2) fails.
\end{example}
This motivates our use of Mémin's theorem below, which allows one to pass to a sufficiently nice subsequence, avoiding the pathologies of the above example.
\begin{proof}
    We show the `only if' implication first. Suppose $\{X^{n}\}_{n}$ admits $Y^{n}\in\mathrm{co}\{X^{m}:m\geq n\}$ such that $\{Y^{n}\}_{n}$ converges in the Émery topology to $Y$. By Mémin's theorem (see Théorème II.3, \cite{memin}) there exists a probability $\mathbb{Q}\sim\mathbb{P}$ and a subsequence $\{n_{k}\}_{k}$ (which we may take to satisfy $n_{k}\geq k$) such that $Y^{n_{k}}=Y^{n_{k}}_{0}+M^{k}+A^{k}$ and $Y=Y_{0}+\widetilde{M}+\widetilde{A}$ where $\{M^{k}\}_{k}$ and $\widetilde{M}$ are $L^{2}(\mathbb{Q})$-martingales, $\{A^{k}\}_{k}$ and $\widetilde{A}$ are predictable with $\mathbb{Q}$-integrable variation, and
    \begin{equation*}
        \lim_{k}\int_{\Omega}\mathrm{var}\left(A^{k}-\widetilde{A}\right)_{T}d\mathbb{Q}=0,
    \end{equation*}
    \begin{equation*}
        \lim_{k}\int_{\Omega}\left[M^{k}-\widetilde{M},M^{k}-\widetilde{M}\right]_{T}d\mathbb{Q}=0.
    \end{equation*}
    Define $\{Z^{k}\}_{k}$ by setting $Z^{k}=Y^{n_{k}}$ for each $k$. Since $n_{k}\geq k$, $Z^{k}\in\mathrm{co}\{X^{j}:j\geq k\}$.
    \par
    We now show that the sequence $\{Z^{n}\}_{n}$ satisfies $(1)\wedge(2)$. Since $\sup_{n}\vert{Z^{n}_{0}}\vert<\infty$, (1) is equivalent to the following. For each $\varepsilon>0$, there exists $K>0$ such that
    \begin{equation*}
        \mathbb{Q}\left(\left\{\sum_{i}\lambda_{i}\left\vert{(\xi^{i}\cdot Z^{i})_{T}}\right\vert>K\right\}\right)<\varepsilon,
    \end{equation*}
    for each $\{\xi_{i}\}_{i}\subset\mathscr{P}(1)$ and $\lambda\in\Delta_{\infty}$. For any $K>0$, Lemma \ref{lem:easy} and Markov's inequality implies that
    \begin{equation*}
        \mathbb{Q}\left(\left\{\sum_{i}\lambda_{i}\left\vert{(\xi^{i}\cdot Z^{i})_{T}}\right\vert>K\right\}\right)\leq\mathbb{Q}\left(\left\{\sum_{i}\lambda_{i}\left\vert{(\xi^{i}\cdot M^{i})_{T}}\right\vert>\frac{K}{2}\right\}\right)
    \end{equation*}
    \begin{equation*}
        +\mathbb{Q}\left(\left\{\sum_{i}\lambda_{i}\left\vert{(\xi^{i}\cdot A^{i})_{T}}\right\vert>\frac{K}{2}\right\}\right)\leq\frac{2}{K}\left(\sum_{i}\lambda_{i}\int_{\Omega}\left(\left\vert{(\xi^{i}\cdot M^{i})_{T}}\right\vert+\left\vert{(\xi^{i}\cdot A^{i})_{T}}\right\vert\right)d\mathbb{Q}\right)
    \end{equation*}
    \begin{equation*}
        \leq\frac{2}{K}\left(\sup_{n}\left(\int_{\Omega}[M^{n},M^{n}]_{T}d\mathbb{Q}\right)^{\frac{1}{2}}+\sup_{n}\int_{\Omega}\mathrm{var}(A^{n})d\mathbb{Q}\right).
    \end{equation*}
    Letting $K$ large enough so
    \begin{equation*}
        \frac{2}{K}\left(\sup_{n}\left(\int_{\Omega}[M^{n},M^{n}]_{T}d\mathbb{Q}\right)^{\frac{1}{2}}+\sup_{n}\int_{\Omega}\mathrm{var}(A^{n})d\mathbb{Q}\right)<\varepsilon
    \end{equation*}
    yields (1). We must now show (2). Let $\{G_{n}\}_{n}$ be a disjoint sequence of predictable sets. It suffices to show that, for each $\varepsilon>0$, there exists $n\in\mathbb{N}$ such that $m\geq n$ implies
    \begin{equation*}
        \sup_{\xi\in\mathscr{P}(1)}\sup_{Z\in\mathrm{co}\{Z^{r}\}_{r}}\mathbb{Q}\left(\left\{{(\xi\mathbf{1}_{G_{m}}\cdot Z)^{\ast}_{T}}>\varepsilon\right\}\right)<\varepsilon,
    \end{equation*}
    Fix $Z=\sum_{i}\lambda_{i}Z^{i}\in\mathrm{co}\{Z^{r}\}_{r}$ (where $\lambda\in\Delta_{\infty}$), $\xi\in\mathscr{P}(1)$, and $m\in\mathbb{N}$. Let $M=\sum_{i}\lambda_{i}M^{i}$, $A=\sum_{i}\lambda_{i}A^{i}$; obviously, $Z=Z_{0}+M+A$. It follows from Markov's inequality, Lemma \ref{lem:easy}, and the Burkholder-Davis-Gundy inequality that
    \begin{equation*}
        \mathbb{Q}\left(\left\{{(\xi\mathbf{1}_{G_{m}}\cdot Z)^{\ast}_{T}}>\varepsilon\right\}\right)\leq\mathbb{Q}\left(\left\{{(\xi\mathbf{1}_{G_{m}}\cdot M)^{\ast}_{T}}>\frac{\varepsilon}{2}\right\}\right)+\mathbb{Q}\left(\left\{{(\xi\mathbf{1}_{G_{m}}\cdot A)^{\ast}_{T}}>\frac{\varepsilon}{2}\right\}\right)
    \end{equation*}
    \begin{equation*}
        \leq\frac{2}{\varepsilon}\left(\int_{\Omega}\left({(\xi\mathbf{1}_{G_{m}}\cdot M)_{T}^{\ast}}+{(\xi\mathbf{1}_{G_{m}}\cdot A)_{T}^{\ast}}\right)d\mathbb{Q}\right)
    \end{equation*}
    \begin{equation*}
        \leq\frac{2}{\varepsilon}\left(C\left(\int_{\Omega}\int_{0}^{T}\mathbf{1}_{G_{m}}d[M,M]d\mathbb{Q}\right)^{\frac{1}{2}}+\int_{\Omega}(\mathbf{1}_{G_{m}}\cdot \mathrm{var}(A))_{T}d\mathbb{Q}\right)
    \end{equation*}
    for some constant $C>0$. Thus it suffices to show $(\mathrm{i})\wedge(\mathrm{ii})$, where:
    \begin{enumerate}[i.]
        \item $\lim_{n}\int_{\Omega}(\mathbf{1}_{G_{n}}\cdot\mathrm{var}(A))_{T}d\mathbb{Q}=0$ uniformly in $A\in\mathrm{co}\{A^{r}\}_{r}$.
        \item $\lim_{n}\int_{\Omega}\int_{0}^{T}\mathbf{1}_{G_{n}}d[M,M]d\mathbb{Q}=0$ uniformly in $M\in\mathrm{co}\{M^{r}\}_{r}$.
    \end{enumerate}
    By the Vitali-Hahn-Saks theorem applied to the sequence of Doléans measures obtainable from $\{\mathrm{var}(A^{n})\}_{n}$, (i) holds after applying the triangle inequality; thus, it suffices to show (ii).
    \par
    Fix an element $M$ of $\mathrm{co}\{M^{r}\}_{r}$. We can write $M=\sum_{i}\lambda_{i}M^{i}$ for some $\lambda\in\Delta_{\infty}$. Let $m\in\mathbb{N}$ be such that $\lambda_{i}=0$ for all $i>m$. Let $\{u_{i}\}_{i}$ and $\{v_{j}\}_{j}$ be $m$-tuples in a Hilbert space $H$ with inner product $\langle{\cdot,\cdot}\rangle$. We have
    \begin{equation}\label{eq:grothtypeineq}
        \left\vert\sum_{i}\sum_{j}\lambda_{i}\lambda_{j}\langle{u_{i},v_{j}}\rangle\right\vert\leq\sup_{i}\Vert{u_{i}}\Vert_{H}\sup_{j}\Vert{v_{j}}\Vert_{H},
    \end{equation}
     from the Cauchy-Schwarz inequality.
     \par
     Straightforward calculation yields
     \begin{equation*}
     \int_{\Omega}\int_{0}^{T}\mathbf{1}_{G_{n}}d[M,M]d\mathbb{Q}=\int_{\Omega}\left[\sum_{i}\lambda_{i}(\mathbf{1}_{G_{n}}\cdot M^{i}),\sum_{j}\lambda_{j}(\mathbf{1}_{G_{n}}\cdot M^{j})\right]_{T}d\mathbb{Q}
     \end{equation*}
     \begin{equation*}
         =\sum_{i}\sum_{j}\lambda_{i}\lambda_{j}\int_{\Omega}\left[(\mathbf{1}_{G_{n}}\cdot M^{i}),(\mathbf{1}_{G_{n}}\cdot M^{j})\right]_{T}d\mathbb{Q}
     \end{equation*}
     \begin{equation*}
         \leq\sup_{i}\int_{\Omega}\left[(\mathbf{1}_{G_{n}}\cdot M^{i}),(\mathbf{1}_{G_{n}}\cdot M^{i})\right]_{T}d\mathbb{Q}
     \end{equation*}
     where the last inequality follows from (\ref{eq:grothtypeineq}). Thus, (ii) is equivalent to (iii), where:
     \begin{enumerate}[iii.]
        \item $\lim_{n}\int_{\Omega}\int_{0}^{T}\mathbf{1}_{G_{n}}d[M,M]d\mathbb{Q}=0$ uniformly in $M\in\{M^{r}\}_{r}$.
     \end{enumerate}
     By the Vitali-Hahn-Saks theorem applied to the sequence of Doléans measures obtainable from $\{[M^{n},M^{n}]\}_{n}$, (iii) holds. This concludes the proof of necessity.
     \par
     We now show sufficiency. Let $Z^{n}\in\mathrm{co}\{X^{m}:m\geq n\}$ be a sequence satisfying $(1)\wedge (2)$. By Theorem \ref{thm:amazedisjoint} and the validity of (1), there exists a disjoint sequence $\{G_{n}\}_{n}$ of predictable sets, and $Y^{n}\in\mathrm{co}\{Z^{m}:m\geq n\}\subset\mathrm{co}\{X^{m}:m\geq n\}$ such that $\{(\mathbf{1}_{(\Omega\times[0,T])\setminus G_{n}}\cdot Y^{n})\}_{n}$ converges in the semimartingale topology to a semimartingale $\widetilde{Y}$. It is no loss of generality to suppose that $\{Y^{n}_{0}\}_{n}$ converges to some $Y_{0}$; define $Y=Y_{0}+\widetilde{Y}$. We claim that $\{Y^{n}\}_{n}$ converges to $Y$ in the semimartingale topology, which would prove the claim.
     \par
     By assumption (2),
     \begin{equation*}
         \lim_{n}\sup_{Z\in\mathrm{co}\{Z^{m}\}_{m}}\mathbf{D}((\mathbf{1}_{G_{n}}\cdot Z),0)=0,
     \end{equation*}
     and in particular,
     \begin{equation*}
         \lim_{n}\sup_{m}\mathbf{D}((\mathbf{1}_{G_{n}}\cdot Y^{m}),0)=0.
     \end{equation*}
     For each $\varepsilon>0$, we therefore may find $m_{\varepsilon}$ such that $\ell\geq m_{\varepsilon}$ implies
     \begin{equation*}
         \mathbf{D}((\mathbf{1}_{G_{\ell}}\cdot Y^{n}),0)<\varepsilon,
     \end{equation*}
     for every $n$. For every $\varepsilon>0$, we may find $w_{\varepsilon}$ such that $\ell\geq w_{\varepsilon}$ implies
     \begin{equation*}
         \mathbf{D}\left(Y^{\ell}_{0}+(\mathbf{1}_{(\Omega\times[0,T])\setminus G_{\ell}}\cdot Y^{\ell}),Y\right)<\varepsilon.
     \end{equation*}
     Let $\ell$ be arbitrary. By the triangle inequality and translation invariance of $\mathbf{D}$,
     \begin{equation*}
         \mathbf{D}(Y,Y^{\ell})\leq\mathbf{D}(Y^{\ell}_{0}+(\mathbf{1}_{(\Omega\times[0,T])\setminus G_{\ell}}\cdot Y^{\ell}),Y)+\mathbf{D}(Y^{\ell},Y^{\ell}_{0}+(\mathbf{1}_{(\Omega\times[0,T])\setminus G_{\ell}}\cdot Y^{\ell}))
     \end{equation*}
     \begin{equation*}
         =\mathbf{D}(Y^{\ell}_{0}+(\mathbf{1}_{(\Omega\times[0,T])\setminus G_{\ell}}\cdot Y^{\ell}),Y)+\mathbf{D}((\mathbf{1}_{\Omega\times[0,T]}\cdot Y^{\ell}),(\mathbf{1}_{(\Omega\times[0,T])\setminus G_{\ell}}\cdot Y^{\ell}))
     \end{equation*}
     \begin{equation}\label{eq:dinequality}
         =\mathbf{D}(Y^{\ell}_{0}+(\mathbf{1}_{(\Omega\times[0,T])\setminus G_{\ell}}\cdot Y^{\ell}),Y)+\mathbf{D}((\mathbf{1}_{G_{\ell}}\cdot Y^{\ell}),0).
     \end{equation}
     Let $\varepsilon>0$ be arbitrary. If $\ell\geq w_{\varepsilon/2}\vee m_{\varepsilon/2}$, then
     \begin{equation*}
         \mathbf{D}(Y,Y^{\ell})<\varepsilon,
     \end{equation*}
     by (\ref{eq:dinequality}) and the definition of $w$ and $m$, which proves the claim.
\end{proof}
\section{Convex predictable uniform tightness}\label{sec:super}
The validity of Theorem \ref{thm:amaze} and Theorem \ref{thm:amazedisjoint} for $\{X^{n}\}_{n}$ is ensured whenever
\begin{equation*}
    \mathrm{co}\left\{\vert{(\xi\cdot X^{n})_{T}\vert}:n\in\mathbb{N},\xi\in\mathscr{P}(1)\right\},
\end{equation*}
is bounded in probability. This condition is a convexified version of \textit{predictable uniform tightness}, which asks that
\begin{equation*}
    \left\{(\xi\cdot X^{n})_{T}:n\in\mathbb{N},\xi\in\mathscr{P}(1)\right\},
\end{equation*}
is bounded in probability. Predictable uniform tightness is a relatively weak condition. In this section, we study whether this relative weakness passes to the convexified version.
\par
Our study is motivated by the failure of Theorem \ref{thm:amaze} and Theorem \ref{thm:amazedisjoint} if one only assumes predictable uniform tightness of $\{X^{n}\}_{n}$.
\begin{example}\label{ex:stable}
    Suppose $T=1$, $(\Omega,\mathscr{F}_{1/2},\mathbb{P})$ is a Lebesgue-Rokhlin probability space, and $\mathscr{F}_{t}=\mathscr{F}_{0}$ for $t<\frac{1}{2}$. Let $\{g_{n}\}_{n}$ be an i.i.d. sequence of $1$-stable random variables. Define, for each $n$, $X^{n}=\mathbf{1}_{\llbracket{1/2,1}\rrbracket}g_{n}$. For essentially the same reasons given in the appendix of \cite{ftap}, $\{X^{n}\}_{n}$ cannot admit convex combinations converging in the sense of Theorem \ref{thm:amazedisjoint}, despite satisfying predictable uniform tightness.
\end{example}
The main result of this section is the following, which the author believes is exceptionally surprising. It applies, after passing to a subsequence, whenever a sequence of supermartingales u.c.p. converges.
\begin{theorem}\label{thm:superucp}
    Suppose $\{X^{n}\}_{n}$ is a sequence of supermartingales such that
    \begin{equation*}
        \mathrm{co}\{\left(X^{n}\right)^{\ast}_{T}:n\in\mathbb{N}\},
    \end{equation*}
    is bounded in probability. Then there exists $Y^{n}\in\mathrm{co}\{X^{m}:m\geq n\}$ such that
    \begin{equation*}
    \mathrm{co}\left\{\vert{(\xi\cdot Y^{n})_{T}\vert}:n\in\mathbb{N},\xi\in\mathscr{P}(1)\right\},
\end{equation*}
    is bounded in probability.
\end{theorem}
\begin{remark}
    If $\{(X^{n})^{\ast}_{T}\}_{n}$ is dominated above by some random variable, one does not need to pass to convex combinations in the statement of Theorem \ref{thm:superucp}.
\end{remark}
Non-convexified versions of Theorem \ref{thm:superucp}, such as Lemma 4.7 of \cite{emeryconvcomp}, are consequences of Burkholder's inequality for supermartingales (see Theorem 47, \cite{meyer-mart-int}). Unfortunately, Burkholder's inequality fails to satisfy convexity, so we must adopt an alternative approach.
\par
Considering only supermartingales is justified by the immense applicability of predictable uniform tightness for sequences of supermartingales (see e.g. \cite{emeryconvcomp}) and the necessity, sans additional assumptions, of the supermartingale property of $\{X^{n}\}_{n}$ for the validity of Theorem \ref{thm:superucp}.
\begin{example}
    Suppose the stochastic basis admits a Brownian motion $B$. Let $X^{n}$ be the strong solution to the Itô stochastic differential equation
    \begin{equation*}
        dX^{n}_{t}=-n^{3}X^{n}_{t}dt+ndB_{t},
    \end{equation*}
    \begin{equation*}
        X^{n}_{0}=0.
    \end{equation*}
    By (Theorem 2.5, \cite{ou}), $\{X^{n}\}_{n}$ u.c.p. converges to zero, and so satisfies the conditions of Theorem \ref{thm:superucp} after passing to a subsequence. However, $\{X^{n}\}_{n}$ does not even admit convex combinations $\{Y^{n}\}_{n}$ satisfying predictable uniform tightness; indeed, $\{[Y^{n},Y^{n}]_{T}\}_{n}$ will always be unbounded in probability, contradicting (Proposition A.1, \cite{emeryconvcomp}). Thus, $\{X^{n}\}_{n}$ fails the conclusion of Theorem \ref{thm:superucp}.
\end{example}
We now prove Theorem \ref{thm:superucp}.
\begin{proof}
    Let $W>0$ be a finite upper bound on $\{\vert{X^{n}_{0}\vert}\}_{n}$, which exists by assumption. Using (Lemma 2.3, \cite{bipolar}), find a probability measure $\mathbb{Q}\sim\mathbb{P}$ such that
    \begin{equation*}
        \sup_{n}\int_{\Omega}\left(X^{n}\right)^{\ast}_{T}d\mathbb{Q}<\infty.
    \end{equation*}
    By (Lemma 2.5, \cite{schruessdies}), we therefore may find a random variable $\zeta>0$ and convex combinations $Y^{n}\in\mathrm{co}\{X^{m}:m\geq n\}$ such that
    \begin{equation*}
        \left(Y^{n}\right)^{\ast}_{T}\leq\zeta,
    \end{equation*}
    for all $n$. Define a sequence $\{\tau_{n}\}_{n}$ of stopping times by
    \begin{equation*}
        \tau_{n}=\inf\left\{t\in[0,T]:\sup_{m}\left(Y^{m}\right)^{\ast}_{t-}\geq n\right\}\wedge T.
    \end{equation*}
    Since $\sup_{m}\left(Y^{m}\right)^{\ast}_{t-}$ is a predictable process, each $\tau_{n}$ is a predictable stopping time. Remark that $\mathbb{P}\left(\{\tau_{n}<T\}\right)\leq\mathbb{P}\left(\{\zeta\geq n\}\right)$, so that
    \begin{equation}\label{eq:smallprob}
        \limsup_{n}\mathbb{P}\left(\{\tau_{n}<T\}\right)=0.
    \end{equation}
    Similarly, for each $\varepsilon>0$, we may find $\widetilde{K}_{\varepsilon}>0$ with $\mathbb{P}\left(\left\{\zeta>\widetilde{K}_{\varepsilon}\right\}\right)<\varepsilon$.
    \par
    It suffices to show the following. For each $\varepsilon>0$, there exists $K>0$ such that
    \begin{equation*}
        \mathbb{P}\left(\left\{\sum_{i}\lambda_{i}\left\vert{(\xi^{i}\cdot Y^{i})_{T}}\right\vert>K\right\}\right)<\varepsilon,
    \end{equation*}
    for each $\{\xi_{i}\}_{i}\subset\mathscr{P}(1)$ and $\lambda\in\Delta_{\infty}$. Fix $\varepsilon>0$. Use (\ref{eq:smallprob}) to find $n$ with $\mathbb{P}(\{\tau_{n}<T\})<\frac{\varepsilon}{3}$. For any $K>0$, Lemma \ref{lem:easy} implies that
    \begin{equation*}
        \mathbb{P}\left(\left\{\sum_{i}\lambda_{i}\left\vert{(\xi^{i}\cdot Y^{i})_{T}}\right\vert>K\right\}\right)\leq\mathbb{P}(\{\tau_{n}<T\})
    \end{equation*}
    \begin{equation*}
        +\mathbb{P}\left(\{\tau_{n}=T\}\cap\left\{\sum_{i}\lambda_{i}\left\vert{(\xi^{i}\cdot Y^{i})_{T}}\right\vert>K\right\}\right)<\frac{\varepsilon}{3}
    \end{equation*}
    \begin{equation}\label{eq:supermartatstop}
        +\mathbb{P}\left(\left\{\sum_{i}\lambda_{i}\left\vert{(\xi^{i}\cdot Z^{i,n})_{T}}\right\vert>\frac{K}{2}\right\}\right)+\mathbb{P}\left(\left\{\sum_{i}\lambda_{i}\left\vert{\Delta{Y^{i}_{\tau_{n}}}}\right\vert>\frac{K}{2}\right\}\right),
    \end{equation}
    where $Z^{i,n}=(Y^{i})^{\tau_{n}}-\mathbf{1}_{\llbracket{\tau_{n},\infty}\llbracket}\Delta{Y^{i}_{\tau_{n}}}$. We may estimate,
    \begin{equation}\label{eq:linfinitybound}
        \left\Vert{(Z^{i,n})^{\ast}_{T}}\right\Vert_{L^{\infty}}\leq n.
    \end{equation}
     By the predictability of $\tau_{n}$ and (Lemma 2.27, \cite{jacshir13}), $Z^{i,n}$ is a supermartingale under $\mathbb{P}$, and therefore admits a Doob-Meyer decomposition $Z^{i,n}-Z^{i,n}_{0}=M^{i,n}-A^{i,n}$, where $M^{i,n}$ is a uniformly integrable martingale under $\mathbb{P}$, and $A^{i,n}$ is a predictable and increasing process (we may assume both start at zero). By (\ref{eq:linfinitybound}) and Meyer's inequality (see Theorem A.2, \cite{sharpbdg}),
     \begin{equation}\label{eq:l2boundmart}
         \left(\int_{\Omega}\left\vert{M^{i,n}_{T}}\right\vert^{2}d\mathbb{P}\right)^{\frac{1}{2}}\leq\sqrt{18}(n+W).
     \end{equation}
     By (\ref{eq:l2boundmart}) and (\ref{eq:linfinitybound}),
     \begin{equation}\label{eq:l1varboundsuper}
         \int_{\Omega}\mathrm{var}(A^{i,n})_{T}d\mathbb{P}=\int_{\Omega}A^{i,n}_{T}d\mathbb{P}\leq(\sqrt{18}+1)(n+W).
     \end{equation}
     By Lemma \ref{lem:easy}, equations (\ref{eq:l2boundmart}) and (\ref{eq:l1varboundsuper}), and Markov's inequality,
     \begin{equation*}
         \mathbb{P}\left(\left\{\sum_{i}\lambda_{i}\left\vert{(\xi^{i}\cdot Z^{i,n})_{T}}\right\vert>\frac{K}{2}\right\}\right)\leq\mathbb{P}\left(\left\{\sum_{i}\lambda_{i}\left\vert{(\xi^{i}\cdot M^{i,n})_{T}}\right\vert>\frac{K}{4}\right\}\right)
     \end{equation*}
     \begin{equation*}
         +\mathbb{P}\left(\left\{\sum_{i}\lambda_{i}\left\vert{(\xi^{i}\cdot A^{i,n})_{T}}\right\vert>\frac{K}{4}\right\}\right)\leq\frac{4(2\sqrt{18}+1)(n+W)}{K}.
     \end{equation*}
     In light of (\ref{eq:supermartatstop}) and the above, letting $K\geq 4\widetilde{K}_{\varepsilon/3}$ large enough so
     \begin{equation*}
         \frac{4(2\sqrt{18}+1)(n+W)}{K}<\frac{\varepsilon}{3},
     \end{equation*}
     yields the claim.
\end{proof}
\begin{remark}\label{rem:localsuper}
    It is not difficult to see from the above proof that Theorem \ref{thm:superucp} holds also for local supermartingales.
\end{remark}
An interesting corollary to Theorem \ref{thm:superucp} is the following version of Helly's selection theorem. It refines or generalizes results due to \cite{schcamp,martconvcomp,melnikov}.
\begin{corollary*}
    Suppose $\{X^{n}\}_{n}$ is a sequence of supermartingales such that
    \begin{equation*}
        \mathrm{co}\{\left(X^{n}\right)^{\ast}_{T}:n\in\mathbb{N}\},
    \end{equation*}
    is bounded in probability. Then there exists $Y^{n}\in\mathrm{co}\{X^{m}:m\geq n\}$ and an optional process $Y$ such that $\{Y^{n}\}_{n}$ converges to $Y$ pointwise outside of an evanescent set.
\end{corollary*}
The supermartingale property plays a significant role. Indeed, without additional assumptions, the above corollary does not generalize to semimartingales. 
\begin{example}
    Let $\{f_{n}\}_{n}$ be an $\ell_{1}$-basis in the unit sphere of $C[0,T]$ (e.g. take the image of the coordinate basis in $\ell_{1}$ under the isometric embedding $\ell_{1}\hooklongrightarrow C[0,T]$ obtained from the Banach-Mazur theorem). For each $n$ we may find an absolutely continuous $g_{n}\in C[0,T]$ with
    \begin{equation*}
        \sup_{t\in[0,T]}\vert{f_{n}(t)-g_{n}(t)}\vert<\frac{1}{n}.
    \end{equation*}
    It is not difficult to see that $\{g_{n}\}_{n}$ does not admit convex combinations converging pointwise on $[0,T]$. Thus, the sequence $\{X^{n}\}_{n}$ of semimartingales defined by $X^{n}=((t,\omega)\longmapsto g_{n}(t))$ does not satisfy the conclusion of the above corollary, despite the bound $\Vert{(X^{n})^{\ast}_{T}}\Vert_{L^{\infty}}\leq1+\frac{1}{n}\leq 2$.
\end{example}
We now give a proof of the corollary.
\begin{proof}
    In light of Theorem \ref{thm:superucp}, it is no loss of generality to assume that
    \begin{equation*}
    \mathrm{co}\left\{\vert{(\xi\cdot X^{n})_{T}\vert}:n\in\mathbb{N},\xi\in\mathscr{P}(1)\right\}
    \end{equation*}
    is bounded in probability. By Lemma \ref{lem:removerestrict}, we therefore may assume that
    \begin{equation*}
        \mathrm{co}\left\{[X^{n},X^{n}]_{T}:n\in\mathbb{N}\right\},
    \end{equation*}
    is bounded in probability. Thus, we may apply Theorem \ref{thm:technicalmemin} and the Borel-Cantelli lemma to pass to convex combinations (still denoted $\{X^{n}\}_{n}$) which decompose as $X^{n}=X^{n}_{0}+M^{n}+A^{n}$, where $\{(M-M^{n})^{\ast}\}_{n}$ converges to zero $\mathbb{P}$-a.s. for some $M$, $\{X^{n}_{0}\}_{n}$ converges, and $\sup_{n}\int_{\Omega}\mathrm{var}(A^{n})_{T}d\mathbb{Q}<\infty$ for some $\mathbb{Q}\sim\mathbb{P}$. Applying (Proposition 13, \cite{schcamp}) and passing to convex combinations yields the claim (c.f. Remark \ref{rem:singularpart}).
\end{proof}
\section*{Acknowledgment}
The author would like to acknowledge and thank the anonymous referee, whose cogent suggestions significantly improved the article.
\printbibliography

@article{kardzit,
    author = {Constantinos Kardaras and Gordan Žitković},
    title = {Forward-convex convergence in probability of sequences of nonnegative random variables},
    journal = {Proceedings of the American Mathematical Society},
    year = {2013},
    volume = {141},
    issue = {3},
    pages = {919-929}
}

@article{sigmaloc,
    author = {Jan Kallsen},
    title = {$\sigma$-localization and $\sigma$-martingales},
    journal = {Theory of Probability and Its Applications},
    year = {2004},
    volume = {48},
    pages = {152-163}
}

@book{ks91,
    author = {Ioannis Karatzas and Steven Shreve},
    title = {Brownian motion and stochastic calculus},
    publisher = {Springer},
    year = {1991},
    address = {New York}
}

@book{jacshir13,
    author = {Jean Jacod and Albert Shiryaev},
    title = {Limit theorems for stochastic processes},
    publisher = {Springer},
    year = {2003},
    address = {Berlin, Heidelberg}
}

@InProceedings{del-sch,
    author = {Freddy Delbaen and Walter Schachermayer},
    title = {A compactness principle for bounded sequences of martingales with applications},
    booktitle = {Seminar on Stochastic Analysis, Random Fields and Applications},
    year = {1999},
    publisher = {Birkhäuser},
    pages = {137-173},
    address = {Basel}
}

@incollection{emery,
    author = {Michel Émery},
    title = {Une topologie sur l'espace des semimartingales},
    booktitle = {Séminaire de Probabilités XIII},
    publisher = {Springer},
    year = {1979},
    pages = {260-280},
    address = {Berlin, Heidelberg}
}

@article{martconvcomp,
    author = {Christoph Czichowsky and Walter Schachermayer},
    title = {Strong supermartingales and limits of nonnegative martingales},
    journal = {Annals of Probability},
    year = {2016},
    volume = {44},
    pages = {171-205}
}

@incollection{bipolar,
    author = {Werner Brannath and Walter Schachermayer},
    title = {A bipolar theorem for ${L}^{0}_{+}({\Omega},\mathcal{F},\mathbf{P})$},
    booktitle = {Séminaire de Probabilités XXXIII},
    publisher = {Springer},
    year = {1999},
    pages = {349-354},
    address = {Berlin, Heidelberg}
}

@article{doobmey,
    author = {Mathias Beiglböck and Walter Schachermayer and Bezirgen Veliyev},
    title = {A short proof of the {Doob-Meyer} theorem},
    journal = {Stochastic Processes and their Applications},
    year = {2012},
    pages = {1204-1209},
    volume = {122}
}

@article{arb,
    author = {Mathias Beiglböck and Walter Schachermayer and Bezirgen Veliyev},
    title = {A direct proof of the {Bichteler-Dellacherie} theorem and connections to arbitrage},
    journal = {Annals of Probability},
    year = {2011},
    pages = {2424-2440},
    volume = {39}
}

@article{ftap,
    author = {Freddy Delbaen and Walter Schachermayer},
    title = {A general version of the fundamental theorem of asset pricing},
    journal = {Mathematische Annalen},
    year = {1994},
    volume = {300},
    pages = {463-520}
}

@article{pure,
    author = {Aleš Černý and Johannes Ruf},
    title = {Pure-jump semimartingales},
    journal = {Bernoulli},
    year = {2021},
    volume = {27},
    pages = {2624-2648}
}

@article{memin,
    author = {Jean Mémin},
    title = {Espaces de semi martingales et changement de probabilité},
    journal = {Zeitschrift für Wahrscheinlichkeitstheorie und Verwandte Gebiete},
    year = {1980},
    volume = {52},
    pages = {9-39}
}

@article{ou,
    author = {Svend Graversen and Goran Peskir},
    title = {Maximal inequalities for the {Ornstein-Uhlenbeck} process},
    journal = {Proceedings of the American Mathematical Society},
    year = {2000},
    volume = {128},
    pages = {3035-3041}
}

@article{schcamp,
    author = {Luciano Campi and Walter Schachermayer},
    title = {A super-replication theorem in {K}abanov's model of transaction costs},
    journal = {Finance and Stochastics},
    year = {2006},
    pages = {579-596},
    volume = {10},
    issue = {4}
}

@article{schruessdies,
    author = {Joseph Diestel and Wolfgang Ruess and Walter Schachermayer},
    title = {Weak compactness in {$L^{1}(\mu,X)$}},
    journal = {Proceedings of the American Mathematical Society},
    year = {1993},
    volume = {118},
    number = {2},
    pages = {447-453},
}

@article{kardem,
    author = {Constantinos Kardaras},
    title = {On the closure in the {Emery} topology of semimartingale wealth-process sets},
    journal = {Annals of Applied Probability},
    year = {2013},
    volume = {23},
    pages = {1355-1376}
}

@article{optdecompconstrain,
    author = {Hans Föllmer and Dmitry Kramkov},
    title = {Optional decompositions under constraints},
    journal = {Probability Theory and Related Fields},
    year = {1997},
    pages = {1-25},
    volume = {109}
}

@article{emeryconvcomp,
    author = {Christa Cuchiero and Josef Teichmann},
    title = {A convergence result for the {Emery} topology and a variant of the proof of the fundamental theorem of asset pricing},
    journal = {Finance and Stochastics},
    year = {2015},
    pages = {743–761},
    volume = {19}
}

@incollection{embound,
    author = {Christophe Stricker},
    title = {Lois de semimartingales et critères de compacité},
    booktitle = {Séminaire de Probabilités XIX},
    publisher = {Springer},
    year = {1985},
    pages = {209-217},
    address = {Berlin, Heidelberg}
}

@article{ogkp,
    author = {Mikhail Kadec and Alexander Pełczyński},
    title = {Bases, lacunary sequences and complemented subspaces in the spaces {$L_{p}$}},
    journal = {Studia Mathematica},
    year = {1962},
    volume = {21},
    pages = {161-176}
}

@article{predictablernderivative,
    author = {Freddy Delbaen and Walter Schachermayer},
    title = {The existence of absolutely continuous local martingale measures},
    journal = {Annals of Applied Probability},
    year = {1995},
    volume = {5},
    pages = {926-945}
}

@book{dellmey,
    author = {Claude Dellacherie and Paul-André Meyer},
    title = {Probabilities and potential},
    publisher = {North-Holland},
    year = {1978},
    address = {Amsterdam, New York, Oxford}
}

@incollection{constraintappl,
    author = {Christoph Czichowsky and Nicholas Westray and Harry Zheng},
    title = {Convergence in the semimartingale topology and constrained portfolios},
    booktitle = {Séminaire de Probabilités XLIII},
    publisher = {Springer},
    year = {2011},
    pages = {395–412},
    address = {Berlin, Heidelberg}
}

@incollection{constraintapplv2,
    author = {Christoph Czichowsky and Martin Schweizer},
    title = {Closedness in the semimartingale topology for spaces of stochastic integrals with constrained integrands},
    booktitle = {Séminaire de Probabilités XLIII},
    publisher = {Springer},
    year = {2011},
    pages = {413–436},
    address = {Berlin, Heidelberg}
}

@article{karatzit,
    author = {Ioannis Karatzas and Gordan Žitković},
    title = {Optimal consumption from investment and random endowment in incomplete semimartingale markets},
    journal = {Annals of Probability},
    year = {2003},
    volume = {31},
    pages = {1821-1858}
}

@article{kramsch,
    author = {Dmitry Kramkov and Walter Schachermayer},
    title = {The asymptotic elasticity of utility functions and optimal investment in incomplete markets},
    journal = {Annals of Applied Probability},
    year = {1999},
    volume = {9},
    pages = {904-950}
}

@article{ptrfskoro,
    author = {Adam Jakubowski and Jean Mémin and Gilles Pages },
    title = {Convergence en loi des suites d'intégrales stochastiques sur l'espace 
 {$\mathbb{D}^{1}$} de {Skorokhod}},
    journal = {Probability Theory and Related Fields},
    year = {1989},
    volume = {81},
    pages = {111–137}
}

@book{brezis,
    author = {Haïm Brezis},
    title = {Functional analysis, {Sobolev} spaces and partial differential equations},
    publisher = {Springer},
    year = {2011},
    address = {New York}
}

@article{unboundftap,
    author = {Freddy Delbaen and Walter Schachermayer},
    title = {The fundamental theorem of asset pricing for unbounded stochastic processes},
    journal = {Mathematische Annalen},
    year = {1998},
    volume = {312},
    pages = {215-250}
}

@incollection{delmeyyor,
    author = {Claude Dellacherie and Paul-André Meyer and Marc Yor},
    title = {Sur certaines propriétés des espaces de {B}anach {$H^{1}$} et {$BMO$}},
    booktitle = {Séminaire de Probabilités XII},
    publisher = {Springer},
    year = {1978},
    pages = {98-113},
    address = {Berlin, Heidelberg}
}

@article{gaoexample,
    author = {Niushan Gao and Denny Leung and Foivos Xanthos},
    title = {On local convexity in {$\mathbb{L}^{0}$} and switching probability measures},
    journal = {Positivity},
    year = {2023},
    volume = {27},
    pages = {1-27}
}

@book{albkalton,
    author = {Fernando Albiac and Nigel Kalton},
    title = {Topics in {Banach} space theory},
    publisher = {Springer},
    year = {2006},
    address = {New York}
}

@book{meyer-mart-int,
    author = {Paul-André Meyer},
    title = {Martingales and stochastic integrals {I}},
    publisher = {Springer},
    year = {1972},
    address = {Berlin, Heidelberg}
}

@article{sharpbdg,
    author = {Walter Schachermayer and Florian Stebegg},
    title = {The sharp constant for the {Burkholder-Davis-Gundy} inequality and non-smooth pasting},
    journal = {Bernoulli},
    year = {2018},
    volume = {24}
}

@article{melnikov,
   title={Relative weak compactness in infinite-dimensional {Fefferman-Meyer} duality},
   author={Vasily Melnikov},
   journal={Journal of Mathematical Analysis and Applications},
   year={2025},
   volume = {543},
   issue = {1}
}
\appendix
\section{Two counterexamples}\label{app:counter}
The convergence described by Theorem \ref{thm:amaze} or Theorem \ref{thm:amazedisjoint} is weaker than semimartingale convergence. We give two counterexamples that demonstrate this claim.
\begin{example}\label{ex:ez}
    Define $X^{n}=\int_{0}^{\cdot}((1/n)\wedge T)^{-1}\mathbf{1}_{\llbracket{0,(1/n)\wedge T}\rrbracket}ds$ ($ds$ is the Lebesgue measure). $\{X^{n}\}_{n}$ satisfies the boundedness conditions of Theorem \ref{thm:amaze} and Theorem \ref{thm:amazedisjoint}, but does not admit convex combinations converging in the semimartingale topology.
\end{example}
\begin{example}\label{ex:harder}
    Suppose $T=1$. We adapt an example of Delbaen and Schachermayer \cite{del-sch}. Let $\bigcup_{n=1}^{\infty}\{\varepsilon_{n,k}\}_{k=1}^{2^{n-1}}$ be a collection of independent random variables with
    \begin{equation*}
        \mathbb{P}\left(\left\{\varepsilon_{n,k}=-2^{-n}\right\}\right)=1-4^{-n},
    \end{equation*}
    \begin{equation*}
        \mathbb{P}\left(\left\{\varepsilon_{n,k}=2^{n}(1-4^{-n})\right\}\right)=4^{-n}.
    \end{equation*}
    Defining $t_{n,k}=\frac{2k-1}{2^{n}}$ (where $n,k\in\mathbb{N}$, $k\leq 2^{n-1}$), let $M$ be the process with
    \begin{equation*}
        M_{t}=\sum_{(n,k):t_{n,k}\leq t}8^{-n}\varepsilon_{n,k}.
    \end{equation*}
    Let the filtration be the stochastic basis generated by $M$. Define $M$-integrable predictable processes $\{H^{n}\}_{n}$ by
    \begin{equation*}
        H^{n}=\sum_{k=1}^{2^{n-1}}8^{n}\mathbf{1}_{\llbracket{t_{n,k}}\rrbracket}.
    \end{equation*}
    If $X^{n}=(H^{n}\cdot M)$, then $\{X^{n}\}_{n}$ satisfies the boundedness conditions of Theorem \ref{thm:amaze} and Theorem \ref{thm:amazedisjoint}, but does not admit convex combinations converging in the semimartingale topology. Indeed, suppose a sequence $\{Y^{n}\}_{n}$ of convex combinations of $\{X^{n}\}_{n}$ converges to some $X$ in the semimartingale topology. The Borel-Cantelli lemma implies that
    \begin{equation*}
        \mathbb{P}\left(\left\{\forall t\in[0,T],X_{t}=-\frac{t}{2}\right\}\right)=1
    \end{equation*}
    so that $X$ must be a nonzero decreasing process. On the other hand, Mémin's theorem implies that $X$ is a $\sigma$-martingale, a contradiction.
\end{example}
Under suitable boundedness conditions, Theorem \ref{thm:technicalmemin} and equation (\ref{eq:decomppredictcanon}) from Section \ref{sec:prelim} implies every sequence $\{X^{n}\}_{n}$ which does not admit convex combinations converging in the semimartingale topology can be split (after passing to convex combinations and switching the measure) into a convergent martingale part, a potentially non-convergent continuous finite variation part, and a potentially non-convergent predictable finite variation pure jump part. The above examples are therefore universal in the following sense. Example \ref{ex:ez} is such that $X^{n}=(X^{n})^{c}$, while Example \ref{ex:harder} is such that $X^{n}\in\mathscr{V}^{d}$, so that every potential counterexample to convex compactness in the semimartingale topology is essentially an interpolation between Example \ref{ex:ez} and Example \ref{ex:harder}.
\end{document}